\documentclass[12pt]{amsart}
\usepackage{amssymb,amscd,graphics,verbatim,mathrsfs,latexsym,amsmath,amsthm,eucal}
\usepackage{a4wide}
\usepackage{verbatim}

\usepackage{color}

\newcommand{\cun}{\mathcal{C}^\infty}
\newcommand{\cz}{\mathbb{C}}
\newcommand{\nz}{\mathbb{N}}
\newcommand{\rz}{\mathbb{R}}
\newcommand{\zz}{\mathbb{Z}}

\newcommand{\cI}{\mathcal{I}}
\newcommand{\cN}{\mathcal{N}}
\newcommand{\cR}{\mathcal{R}}
\newcommand{\ctm}{{}^cTM}
\newcommand{\ctsm}{{}^cT^*M}

\newcommand{\Diff}{\mathrm{Diff}}
\newcommand{\Dh}{D_{(\partial M,h)}}
\newcommand{\Dom}{\mathcal{D}}
\newcommand{\Dp}{D_\partial}
\newcommand{\End}{\mathrm{End}}
\newcommand{\ind}{\mathrm{index}}
\newcommand{\Mci}{M^\circ}
\newcommand{\Me}{M_\epsilon}

\newcommand{\res}{\mathrm{res}}
\newcommand{\Res}{\mathrm{Res}}
\newcommand{\sus}{\mathrm{sus}}
\newcommand{\tD}{\tilde{D}}

\newcommand{\tU}{\tilde{U}}
\newcommand{\tV}{\tilde{V}}
\newcommand{\Spec}{\mathrm{Spec}}
\newcommand{\tr}{\mathrm{tr}}
\newcommand{\Tr}{\mathrm{Tr}}
\newtheorem{corollary}{Corollary}%[section]
\newtheorem{theorem}[corollary]{Theorem}
\newtheorem*{mtheorem}{Main Theorem}
\newtheorem{proposition}[corollary]{Proposition}
\newtheorem{lemma}[corollary]{Lemma}
\newtheorem{assumption}{Assumption}
\theoremstyle{definition}
\newtheorem{definition}[corollary]{Definition}
\theoremstyle{remark}
\newtheorem{remark}[corollary]{Remark}

\def\dbar{d\hspace*{-0.08em}\bar{}\hspace*{0.1em}}

\def\tr{\text{tr}}

\newcommand{\hh}{\mathbb{H}}

\begin{document}

\title[Regularity of the eta function]{Regularity of the eta function on
manifolds with cusps}
\author{Paul Loya}
\address{Department of Mathematics \\Binghamton University\\
Binghamton\\NY 13902\\U.S.A. }
\email{paul@math.binghamton.edu}
\author{Sergiu Moroianu}
\address{Institutul de Matematic\u{a} al Academiei Rom\^{a}\-ne\\
P.O. Box 1-764\\
RO-014700 Bu\-cha\-rest, Romania}
\email{moroianu@alum.mit.edu}
\author{Jinsung Park}
\address{School of Mathematics\\ Korea Institute for Advanced Study\\
Hoegiro 87\\ Dong\-daemun-gu\\ Seoul 130-722\\
Korea}
\email{jinsung@kias.re.kr}
\begin{abstract}
On a spin manifold with conformal cusps, we prove under an
invertibility condition at infinity that the eta function of the
twisted Dirac operator has at most simple poles and is regular at
the origin. For hyperbolic manifolds of finite volume, the eta
function of the Dirac operator twisted by any homogeneous vector
bundle is shown to be entire.
\end{abstract}
\date{\today}
\maketitle

%%%%%%%%%%%%%%%%%%%%%%%%%%%%%%%%%%%%%%%%%%%%%%%%%%%%%%%%%%%%%%%%%%%%%%
\section{Introduction}
%%%%%%%%%%%%%%%%%%%%%%%%%%%%%%%%%%%%%%%%%%%%%%%%%%%%%%%%%%%%%%%%%%%%%%

The eta invariant was first introduced in \cite{aps} as a real number
associated to certain elliptic first-order differential operators
on compact manifolds with boundary, which happened to equal the difference
between the Atiyah-Singer integral and the index
with respect to the Atiyah-Patodi-Singer spectral boundary condition.
During the thirty years since its discovery,
this invariant has risen from the status of ``error term''
to that of a subtle tool, highly efficient in solving
otherwise intractable problems from various fields of mathematics.
Let us mention in this respect its recent application in finding obstructions
for hyperbolic and flat $3$-manifolds to bounding hyperbolic
$4$-manifolds \cite{longreid1}.

For a Hermitian vector bundle $E$ over a closed manifold $M$,
consider an elliptic self-adjoint first-order differential
operator $D:\cun(M,E)\to \cun(M,E)$. Then the $L^2$ spectrum of
$D$ is purely discrete and distributed according to the classical
Weyl law. It follows that the complex function
\[\eta(D,s):=\sum_{\lambda\in\Spec(D)\setminus \{0\}}\lambda |\lambda|^{-s-1}\]
is well-defined (and holomorphic) when $\Re(s)>\dim(M)$. This
function admits a meromorphic extension to the complex plane with
possible simple poles. If $\dim(M)$ is odd, the possible poles are
located at $\dim(M)-1-2\nz$ where $\nz=\{0,1,2,\ldots\}$. If
$\dim(M)$ is even and $D$ is a Dirac operator associated to a
Clifford connection, the eta function is entire \cite{BraGyl}.

It is a byproduct of the index theorem of Atiyah, Patodi and
Singer that when $M$ is a boundary and $D$ is the tangential part
of an elliptic operator as above, the point $s=0$ is always
regular for the eta function. In fact, the eta function is always
regular at the origin for general elliptic pseudodifferential
operators; this was proved using $K$-theory in the spirit of early
index theory by Atiyah, Patodi, and Singer \cite{aps3} in the
odd-dimensional case, and by Gilkey
\cite{Gilkey} for arbitrary dimensions. The \emph{eta invariant}
of $D$ is defined as
\[\eta(D)=\eta(D,0)+\dim\ker(D).\]

Still on a closed manifold another question arises: it is easy to
note that the possible residue of the eta function at the origin
is the integral on $M$ of a well-defined density determined
locally by $D$, called the local eta residue. Is this density
zero? For arbitrary differential operators the answer is negative,
but for Dirac operators it was proved by Bismut and Freed
\cite{bf} that this is indeed the case. The proof uses the
properties of the heat coefficients in terms of Clifford
filtration, along the lines of Bismut's heat equation proof of the
local index formula.

In this note we first revisit the vanishing of the local eta density
from the point of view of conformal invariance. We give a self-contained proof,
using the APS formula, of the fact that the eta invariant of
the spin Dirac operator is insensible to conformal changes.
This important fact belongs to the mathematical
folklore but we could not find a complete proof in the literature.
The existing proofs (e.g.\ \cite[pp. 420-421]{APSII}) tend to
use the index formula of \cite{aps} for metrics which are not of product type
near the boundary, without explaining why one can do so.
In section \ref{sec-conformal} we show how to apply the APS formula
to a true product-type metric in order to prove the conformal invariance of the
eta invariant. We deduce the vanishing of the local eta residue from this
conformal invariance, by interpreting
the variation of the eta invariant in terms of the Wodzicki residue.

Our main results concern the eta function on noncompact spin
manifolds with conformal cusps, in particular on complete
finite-volume hyperbolic spin manifolds. More precisely, let $M$
be a compact manifold with boundary and ${[0,\epsilon)}_x\times
\partial M$ a collar neighborhood of its boundary. The interior
$\Mci$ of $M$ is called a conformally cusp manifold if it is
endowed with a metric $g_p$ which near $x=0$ takes the form
\begin{equation}\label{e:metric-p}
g_p = x^{2p}\left(\frac{dx^2}{x^4}+h\right),
\end{equation}
for some $p>0$, where $h$ is a metric on $\partial M$ independent
of $x$. The main examples are complete hyperbolic manifolds of
finite volume, for which $p=1$ and $h$ is flat. Assume now that
$M$ is spin. Let $E$ be a bundle with connection over $M$ which is
of product-type on the collar, and let $D_p$ denote the associated
twisted Dirac operator on $\Sigma \otimes E$ where $\Sigma$ is the
spinor bundle over $M$. We assume that the spin structure and the
connection on $E$ are ``nontrivial'' (Assumption \ref{a1} in Section
\ref{s:conf-cusp}) in the sense that the twisted Dirac operator $\Dh$ for
the induced spin structure on $(\partial M, h)$ is invertible.
Under this assumption, the twisted Dirac operator $D_p$ is
essentially self-adjoint with discrete spectrum obeying the Weyl
law and the corresponding eta function $\eta(D_p,s)$ has a
meromorphic extension to $\cz$ with possible double poles
\cite{wlom}. For the untwisted Dirac operator on finite-volume
hyperbolic manifolds, it was already noted by B\"ar \cite{bar}
that the spectrum of $D=D_1$ is discrete if and only if the spin
structure is ``nontrivial'' on the cusps in the above sense,
otherwise the continuous spectrum of $D$ is $\rz$. In Appendix
\ref{sec-spec} we prove that the same occurs for conformal cusp
metrics: If $\Dh$ fails to be invertible, then for $p \leq 1$ the
twisted Dirac operator $D_p$ has essential spectrum equal to
$\mathbb{R}$. In the case $p > 1$, $D_p$ fails to be essentially
self-adjoint and although every self-adjoint extension of $D_p$
has discrete spectrum, nothing is known about the meromorphic
properties of the corresponding eta functions. These are the
reasons the ``nontriviality'' assumption plays an important
r\^{o}le in this theory. The main results of this paper are that
under Assumption \ref{a1} the eta function $\eta(D_p,s)$ in fact
has at most simple poles, and is always regular at the origin.
Moreover, the poles disappear for hyperbolic manifolds, thus the
eta function is entire in that case.

\begin{mtheorem}
Let $\Mci$ be an odd-dimensional spin manifold with conformal
cusps, $E$ a twisting bundle of product type on the cusps, and let
$D_p$ be the associated twisted Dirac operator to
\eqref{e:metric-p} on $\Sigma \otimes E$ satisfying Assumption
\ref{a1}. Then
\begin{enumerate}
\item The eta function $\eta(D_p,s)$ of the twisted spin
Dirac operator is regular for $\Re(s)>-2$ and has at most simple
poles at $s\in \{-2,-4,\ldots\}$.
\item {When $p=1$, $\Mci$ is
hyperbolic of finite volume and $E$ is a homogeneous vector
bundle, the eta function $\eta(D,s)$ of the twisted Dirac operator
is entire.}
\end{enumerate}
\end{mtheorem}

In even dimensions the eta function vanishes identically since the
spectrum is symmetric, see Section \ref{sec-conformal}.

A related question may be asked on more complicated metrics at
infinity, like the fibered-cusp metrics arising on
$\mathbb{Q}$-rank one locally symmetric spaces. A similar problem
arose from \cite{jipa}, where we could obtain a meromorphic
extension of the eta function for cofinite quotients of
$\mathrm{PSL}_2(\rz)$ by using the Selberg trace formula. The
methods employed here do not seem to extend easily to such spaces.

We now outline this paper. We begin in Section \ref{sec-conformal}
by proving that on a closed spin manifold the eta invariants are
identical for two Dirac operators associated to conformal metrics.
In Section \ref{sec-res} we review the Guillemin-Wodzicki residue
density and residue trace and derive some of their elementary
properties that we need in the sequel. In Section \ref{sec-vanish}
we give a new proof that on any spin manifold, the local eta
residue of a twisted Dirac operator vanishes. In Sections
\ref{s:conf-cusp} and \ref{sec-hyper} we prove the main theorem
based on Melrose's cusp calculus \cite{meni}, which we review in
Appendix \ref{sec-cuspcalc}.

%%%%%%%%%%%%%%%%%%%%%%%%%%%%%%%%%%%%%%%%%%%%%%%%%%%%%%%%%%%%%%%%%%%%%%
\section{Conformal invariance of eta invariants on closed manifolds}\label{sec-conformal}
%%%%%%%%%%%%%%%%%%%%%%%%%%%%%%%%%%%%%%%%%%%%%%%%%%%%%%%%%%%%%%%%%%%%%%

The eta invariant of the spin Dirac operator and of the odd
signature operator are known to be invariant under conformal
changes of the metric. Since on one hand we need to understand
this fact in depth, and on the other hand we were unable to find a
good reference, we chose to give here a complete proof.

Let $(M,g)$ be a closed spin Riemannian manifold of dimension $n$, $(E,\nabla)$ a Hermitian
vector bundle on $M$ with compatible connection $\nabla$, and $D$
the twisted Dirac operator. Note that when $n$ is even,
the eta invariant reduces to $\dim\ker(D)$ (which is known to be a conformal
invariant, see \eqref{confch}). Indeed, in even dimensions the operator $D$ is odd with respect to
the splitting in positive and negative spinors. Thus the eta function itself
vanishes in these dimensions because the spectrum is symmetric around $0$.
For the untwisted spin Dirac operator, the same vanishing occurs in dimensions
$4k+1$: for $n=8k+1$ the spinor bundle has a real structure (i.e.\ a skew-complex map $C$
with $C^2=1$) which anti-commutes with $D$, while in dimensions $8k+5$
it has a quaternionic structure (i.e.\ a skew-complex map $J$ with $J^2=-1$)
which anti-commutes with $D$ \cite[pp.\ 61, Remark (3)]{aps}.

Let $f\in\cun(M)$ be a real conformal factor, $g':=e^{-2f}g$ a
metric conformal to $g$ and $D'$ the corresponding Dirac operator.

\begin{proposition}\label{vancc}
The eta invariants of $D$ and $D'$ coincide.
\end{proposition}
\begin{proof}
The map of dilation by $e^f$ gives an $\mathrm{SO}(n)$-isomorphism
between the orthonormal frame bundles of $g$ and $g'$. Thus the
principal $\mathrm{Spin}(n)$-bundle (for the fixed spin structure)
corresponding to $g$ and $g'$ are also isomorphic via the lift of
this map. This identifies the spinor bundles for the two metrics;
the Dirac operators are linked by the formula
\begin{align}\label{confch}
g'=e^{-2f}g,&&D'=e^{\frac{n+1}{2}f}De^{-\frac{n-1}{2}f}
\end{align}
(see e.g.\ \cite[Proposition 1]{nistor} for a proof). In
particular, the null-spaces of these two operators have the same
dimension.

Let $\psi:I=[0,1]\to \rz$ be a smooth function which is $0$ for $t<1/3$
and which is identically $1$ for $t>2/3$. Set $f_t:=\psi(t)f$ and define
a metric on $X:=I\times M$ by
\[h=dt^2+e^{-2f_t}g.\]
We denote again by $E$ the pull back of $E$ from the second factor,
together with its connection. Therefore, the curvature tensor of $E$
on $X$ satisfies
\begin{equation}\label{part}
\partial_t\lrcorner R^E=0.
\end{equation}

The metric $h$ is of product type near $\partial X$ and hence
the Atiyah-Patodi-Singer formula can be applied to the (chiral)
twisted Dirac operator $D^+$ on $X$:
\[\ind(D^+)=\int_X \hat{A}(R^h)\mathrm{ch}(R^E)-\frac12\eta(D)+\frac12\eta(D').
\]

On the other hand, this index is also equal to the spectral flow
in the space of Riemannian metrics from $D$ to $D'$; again by \eqref{confch},
there is no spectral flow so the index vanishes.
The proof will be concluded by showing that the top component of
the integrand in the APS formula
vanishes.

From \eqref{part} we deduce $\partial_t\lrcorner \exp(R^E)=0$ so
it is enough to show that $\partial_t\lrcorner \hat{A}(R^h)=0$.
Recall that $\hat{A}$ is a polynomial in the Pontrjagin forms
$\tr(R^h)^{2k}\in\Omega^{4k}(X)$. Also, recall that the Pontrjagin
forms are conformal invariant (they only depend on the Weyl tensor
-- first proved by Chern and Simons \cite{CS}). Let
\[h'=e^{2f_t}h=e^{2f_t}dt^2+g.\]
We claim that $\partial_t\lrcorner \tr((R^{h'})^{2k})=0$ for all
$k$. Indeed, let $\nabla$ be the Levi-Civita connection of $h'$.
For every vector field $V$ on $M$ denote by $\tV$ its pull-back to
$X$, which is orthogonal on the length-$1$ vector field
$T:=e^{-f_t}\partial_t$. Note that
\begin{align*}
[\tV,T]=-\psi(t)V(f)T,&&[\tV,\tU]=\widetilde{[V,U]}.
\end{align*}
We deduce that
\begin{align*}
2\langle\nabla_{\tV} T,\tU\rangle=
&\tV\langle T,\tU\rangle+T\langle\tV,\tU\rangle-\tU\langle\tV,T\rangle\\
&+\langle[\tV,T],\tU\rangle+\langle[\tU,\tV],T\rangle
+\langle[\tU,T],\tV\rangle\\=&0.
\end{align*}
Clearly, since also $\langle\nabla_{\tV} T,T\rangle=0$, we infer
$\nabla_{\tV} T=0$. Directly from the definition of the curvature
this implies that $R^{h'}_{\tV \tU}T=0$.
If we split $TX$ into $TI\oplus TM$, we see from the symmetry
of the curvature tensor that
$R^{h'}_{\tV \tU}$ is a diagonal linear map, while $R^{h'}_{\tV T}$
is off-diagonal. It follows that
$\partial_T \lrcorner(R^{h'})^{2k}$ is an off-diagonal form-valued
endomorphism (since it contains exactly one curvature term involving $T$).
Hence its trace is zero.
\end{proof}

{\begin{remark} The same proof applies as well to the (twisted)
odd signature operator on any orientable manifold, since the
Hirzebruch $L$-form, like the $\hat{A}$-form, is also a polynomial
in the Pontrjagin forms.
\end{remark}}

%%%%%%%%%%%%%%%%%%%%%%%%%%%%%%%%%%%%%%%%%%%%%%%%%%%%%%%%%%%%%%%%%%%%%%

\section{The residue trace and the residue density}\label{sec-res}
%%%%%%%%%%%%%%%%%%%%%%%%%%%%%%%%%%%%%%%%%%%%%%%%%%%%%%%%%%%%%%%%%%%%%%

We review a refined construction of the residue trace. Let $A$ be
a classical pseudodifferential operator $A$ of integer order on a
smooth closed manifold $M$ of dimension $n$, acting on
the sections of a vector bundle $E$. We will later be interested
in twisted spinor bundles over spin manifolds, but the description
of the residue density does not need these assumptions. Let
$\kappa_A(m,m')$ be the Schwartz kernel of $A$, which is a
distributional section in $E\boxtimes (E^*\otimes|\Omega|)$ over
$M\times M$ with singular support contained in the diagonal.
Choose a diffeomorphism
\begin{align}
\Phi:U\to V\subset M\times M, &&\Phi(m,v)=(m,\phi_m(v)) \label{Phidiff}
\end{align}
from a neighborhood of the zero section in $TM$ to a neighborhood of the
diagonal in $M$, extending the canonical identification of $M$ with the diagonal.
Cut-off $\kappa_A$ away from
the diagonal, i.e., multiply it by a function $\psi$ with support in $V$
which is identically $1$ near the diagonal. Fix a connection in $E$,
so that we can identify $E^*_{\phi_m(v)}$ with $E^*_m$
using parallel transport along the curve $t\mapsto \phi_m(tv)$. Then
$\Phi^*(\psi \kappa_A)$
is a compactly-supported distributional section over $TM$ in the bundle
$\End(E)$ pulled back from the base, tensored with the fiberwise density bundle.
This distribution is conormal to the zero section, thus
by definition there exists a classical symbol $a(m,\xi)$
on $T^*M$ (with values at $(m,\xi)$ in $\End(E_m)$)
such that
\[\Phi^*(\psi\kappa_A)(m,v) = \frac{1}{(2\pi)^n} \int_{T^*M/M} e^{i\xi(v)}a(m,\xi)\, \omega^n.\]
Here $\omega$ is the canonical symplectic form on $T^*M$, and
$\int_{T^*M/M}$ means integration along the fibers of $T^*M$. The
result on the right-hand side is an $\End(E)$-valued density in
the \emph{base} variables; however since the vertical tangent
bundle to $TM$ at $(m,v)$ is canonically isomorphic to $T_mM$,
this can be interpreted as a vertical density.

Let $\cR$ be the radial (vertical) vector field in the fibers of $T^*M$.
Let $a_{[-n]}$ denote the component of homogeneity
$-n$ of the classical symbol $a$. Fix a Euclidean metric $g$ in the vector bundle
$T^*M$ (this amounts to choosing a Riemannian metric on $M$),
thus defining a sphere bundle $S^*M$ inside $T^*M$.

\begin{definition} \label{def-res}
The residue density of $A$ is the smooth $\End(E)$-valued density
\[
\res(A) := \frac{1}{(2\pi)^n} \int_{S^*M/M} a_{[-n]}\,
\cR\lrcorner \omega^n.
\]
\end{definition}

At this stage, $\res(A)$ depends on a number of choices: the embedding
$\Phi$, the cut-off $\psi$, the connection in $E$ and the metric $g$.

One way to show that $\res(A)$ is defined {independently} of the
choices involved is through holomorphic families. Let
$(A_s)_{s\in\cz}$ be a holomorphic family of pseudodifferential
operators on $E$ such that $A_s$ is of order $k-s$, where $k$ is
the order of $A$, and $A_0 - A \in \Psi^{-\infty}(M,E)$. Then for
$\Re(s)$ sufficiently large, restricting the Schwartz kernel of
$A_s$ to the diagonal $\Delta$ gives a well-defined and
holomorphic {$\End(E)$}-valued density
\[F(s):=\kappa_{A_s}|_{\Delta}.\]
This density extends to $\cz$ with possible simple poles at $s\in
n + k -\nz$, where $k$ is the order of $A$ and  $\nz =
\{0,1,2,3,\ldots\}$. One natural choice of a holomorphic family is
$A_s = A Q_s$ or $A_s = Q_s A$ where $(Q_s)_{s\in\cz}$ is a
holomorphic family of pseudodifferential operators on $E$ such
that $Q_s$ is of order $-s$ and
$Q_0-\mathrm{Id}\in\Psi^{-\infty}(M,E)$. It is then
straightforward to check that
\begin{equation} \label{ResFA}
\Res_{s=0}F(s) = \res(A) .
\end{equation}
Thus, the residue $\Res_{s=0}F(s)$ is well-defined irrespective of the choice
of $A_s$, and also $\res(A)$ is well-defined
independently of choices.
The residue trace of $A$ is defined by
\begin{equation}\label{trres}
\Tr_R(A) = \Res_{s = 0} \Tr(A_s) = \int_M \tr(\res(A)).
\end{equation}

Armed with this holomorphic family interpretation of $\res(A)$ and $\Tr_R(A)$,
one can deduce without effort various properties of $\res$ and $\Tr_R$.
For example, it follows that $\Tr_R$ vanishes on commutators.
Indeed, given integer order operators $A$ and $B$ and taking any auxiliary
family $Q_s$ as explained above, we can write
\[
\Tr([A,B] Q_s)  = \Tr(C_s) + \Tr([A Q_s, B]),
\]
where $C_s = A B Q_s - A Q_s B$ is a holomorphic family of
operators that is smoothing at $s = 0$.  One can
check that $\Tr[S,T] = 0$ for any pseudodifferential operators $S$
and $T$ with $\mathrm{ord}(S) + \mathrm{ord}(T) < -n$, so for
$\Re(s)$ sufficiently large, the trace $\Tr([A Q_s, B])$ vanishes.
Therefore, by analytic continuation, $\Tr([A Q_s, B])$ vanishes
for all $s \in \cz$; in particular, we have $\Tr([A,B] Q_s)  =
\Tr(C_s)$.   Now $C_0$ is smoothing, and since the residue density
of a smoothing operator is zero, we have $\Res_{s = 0} \Tr(C_s) =
\Tr_R(C_0) = 0$. Thus, $\Tr_R([A,B]) = 0$.

Furthermore, if $S$ is any section in $\End(E)$, using the fact that
\[\kappa_{S A_s}=S \kappa_{A_s}\quad \text{and}\quad \kappa_{A_s S}=
\kappa_{A_s} S,\]
where $A_s$ is a holomorphic family as in the formula \eqref{ResFA},
we have $\res(SA)=S\res(A)$ and $\res(AS) = \res(A)S$. That $\res(SA)=S\res(A)$ also follows directly from Definition \ref{def-res}.  In particular, we have
\begin{equation}\label{trace}
\res(u A)=\res(Au)=u\, \res(A)
\end{equation}
for every function $u\in\cun(M,\cz)$. Finally, observing that
\[
\kappa_{A_s}|_\Delta = \Big(\kappa_{A_s^*}|_\Delta \Big)^*,
\]
taking the residue at $s = 0$ of both sides we obtain $\res(A) = \res(A^*)^*$;
that is, we have $\res(A^*)=\res(A)^*$.

%%%%%%%%%%%%%%%%%%%%%%%%%%%%%%%%%%%%%%%%%%%%%%%%%%%%%%%%%%%%%%%%%%%%%%
\section{Vanishing of the local eta residue}\label{sec-vanish}
%%%%%%%%%%%%%%%%%%%%%%%%%%%%%%%%%%%%%%%%%%%%%%%%%%%%%%%%%%%%%%%%%%%%%%

Let $(D_t)_{t\in I}$ be a smooth $1$-parameter family of elliptic self-adjoint
pseudodifferential operators of order $1$ on a closed manifold $M$.
For simplicity assume for a moment that $D_0$ is invertible, hence $D_t$
is invertible for small enough $t$. Then
\begin{align*}
\partial_t \eta(D_t,s)=&\partial_t\Tr\left(D_t (D_t^2)^{-\frac{s+1}{2}}\right)\\
=&\Tr\left(\dot{D}_t (D_t^2)^{-\frac{s+1}{2}}\right)-\frac{s+1}{2}
\Tr\left(D_t (D_t^2)^{-\frac{s+3}{2}}(\dot{D}_tD_t+D_t\dot{D}_t)\right)\\
=&-s\Tr\left(\dot{D}_t (D_t^2)^{-\frac{s+1}{2}}\right)\\
=&-s\Tr\left(\dot{D}_t |D_t|^{-1}(D_t^2)^{-\frac{s}{2}}\right).
\end{align*}
Now $Q_s:=(D_t^2)^{-\frac{s}{2}}$ is an analytic family of
operators of order $-s$ and $Q_0=\mathrm{Id}$, therefore
\begin{equation}\label{vareta}
\partial_t \eta(D_t)=\left[-s\Tr\left(\dot{D}_t |D_t|^{-1}
Q_s\right)\right]_{s=0} = - \Tr_R(\dot{D}_t |D_t|^{-1}).
\end{equation}
From \eqref{trres}, the Wodzicki residue trace vanishes on
smoothing operators, so as a corollary we see that the eta
invariant is constant under smoothing perturbations. By this
argument, the above expression makes sense even when $D_t$ is not
invertible.

In the same spirit, let $D$ be an elliptic
self-adjoint invertible pseudodifferential operator of order
$k\in(0,\infty)$. Then the residue at $s=0$ of the eta function
$\eta(D,s)$ is
\[\Res_{s=0}\Tr(D|D|^{-1}(D^2)^{-\frac{s}{2}})=
\frac{1}{k}\Tr_R(D|D|^{-1}) =\frac{1}{k}\int_M\tr(\res({D}
|D|^{-1})).\]

We have been
assuming that $M$ is closed, so that the trace on the left is
defined. However, notice that by definition, $\tr(\res(D
|D|^{-1}))$ is a local quantity in the sense that it depends only
on finitely many terms of the local symbol of $D|D|^{-1}$;
moreover, each homogeneous term of $D|D|^{-1}$ is given by a
universal formula in terms of the local symbol of $D$ in any
coordinate patch. Using this universal formula for the $-n$ degree
homogeneous term allows us to define the local eta residue on the
\emph{interior} of any manifold (with or without boundary, compact
or not), even in the case that $|D|^{-1}$ does not exist.

\begin{definition}
Let $M$ be a possibly non-compact manifold, $E\to M$ a vector bundle
and $z\in\cz$.
For any elliptic pseudodifferential operator $D\in\Psi^z(M,E)$, the density
\[\tr(\res({D} |D|^{-1}))\in|\Omega(M)|.\]
is called the \emph{local eta residue} of $D$.
\end{definition}
From the definition of the residue density,
the local eta residue is constant under smoothing perturbations,
so the definition makes sense when $D$ is not invertible, non symmetric.
The local eta residue can be
non-vanishing in general (when $M$ is compact,
its integral always vanishes for self-adjoint
operators of positive order since the eta function
is regular at $s=0$ \cite{Gilkey}). However, for Dirac operators we have:

\begin{theorem}\cite{bf}\label{th4}
Let $(M,g)$ be a spin Riemannian manifold and $E$ a twisting bundle.
Then the local eta residue $\tr(\res(D|D|^{-1}))$ of the twisted Dirac operator
vanishes.
\end{theorem}
\begin{proof} We give here a new, easy proof.
Assume that $M$ is closed; this theorem is a local question, so this case suffices.
Let $f$ be an arbitrary smooth real function on $M$. Define $D_t$ as
the Dirac operator associated to the family of conformal metrics
$e^{-2tf}g$. This operator is an unbounded operator in the
$L^2$ space associated to the measure $e^{-ntf}\mu_g$. To work in the fixed
Hilbert space $L^2(\mu_g)$, conjugate through the unitary transformation
\begin{align*}
L^2(M,\Sigma\otimes E,e^{-ntf}\mu_g)\to L^2(M,\Sigma\otimes
E,\mu_g)&& \phi\mapsto e^{-\frac{ntf}{2}}\phi
\end{align*}
where $\Sigma$ denotes the spinor bundle over $M$. Using \eqref{confch},
$D_t$ conjugates to
\begin{align}
\tD_t=&e^{\frac{tf}{2}} D e^{\frac{tf}{2}}
\text{ acting in $L^2(M,\Sigma\otimes E,\mu_g)$} \nonumber\\
\intertext{and we compute} \label{vard}
\partial_t\tD_t=&\frac{1}2 (f\tD_t+\tD_tf).
\end{align}
Using Proposition \ref{vancc} we have on one hand
\[\partial_t\eta(\tD_t)=\partial_t\eta(D_t)=0.\]
On the other hand, plugging \eqref{vard} at $t=0$
into \eqref{vareta} we write:
\begin{align*}
{-}\partial_t \eta(\tD_t)_{|t=0}=&\Tr_R\left[\frac{1}{2} (f D+Df)
|D|^{-1}\right] \quad \quad \text{since $\tD_0=D$}\\
=&\frac{1}{2} \Tr_R(fD|D|^{-1}+f |D|^{-1} D)
\quad \text{since $\Tr_R$ is a trace}\\
=&\Tr_R(fD |D|^{-1}) \qquad \qquad \text{since $D$ commutes with $|D|$}\\
=&\int_M \tr(\res(fD |D|^{-1})).
\end{align*}
From the definition (see also \eqref{trace}), $\res(fD
|D|^{-1})=f\res(D |D|^{-1})$. Since $f$ was arbitrary, we deduce
$\tr(\res(D |D|^{-1}))=0$ as claimed.

\end{proof}

We will need such a vanishing result for a larger class of first-order symmetric
differential operators:

\begin{corollary} \label{regconj}
Let $D$ be a twisted Dirac operator on a spin manifold $(M,g)$.
For any $u\in\cun(M,{\mathbb{R}})$,
the operator $D_u:=e^{-u}De^u$ is symmetric on $M$ with respect to the measure
$\mu_u:= e^{2u} \mu_g$, and the local eta residue of $D_u$ vanishes.
\end{corollary}
\begin{proof}
It is clear that $D_u$ is formally self-adjoint with respect to
the measure $\mu_u$. As in the proof of Theorem \ref{th4}, to
prove that the local eta residue of $D_u$ vanishes we can assume
that $M$ is compact. Then $|D_u|=e^{-u}|D|e^u$, hence
$D_u|D_u|^{-1}=e^{-u}D|D|^{-1}e^u$. By \eqref{trace},
$\res(D_u|D_u|^{-1})=\res(e^{-u}D|D|^{-1}e^u)=\res(D|D|^{-1})$.
The trace of this last endomorphism-valued density vanishes by
Theorem \ref{th4}.
\end{proof}

%%%%%%%%%%%%%%%%%%%%%%%%%%%%%%%%%%%%%%%%%%%%%%%%%%%%%%%%%%%%%%%%%%%%%%

\section{Eta function on conformally cusp manifolds}\label{s:conf-cusp}
%%%%%%%%%%%%%%%%%%%%%%%%%%%%%%%%%%%%%%%%%%%%%%%%%%%%%%%%%%%%%%%%%%%%%%

We turn now to our main object of study.

Let $\Mci$ be the interior of a compact manifold with boundary $
M$ of dimension $n$. We assume that $M$ is spin with a
fixed spin structure, and that the metric is of conformally cusp
type as in \cite{wlom}. To explain this notion, let $x:
M\to[0,\infty)$ be a boundary-defining function for the smooth
structure of $M$, namely
\begin{enumerate}
\item $x\in\cun( M)$; \item $\{x=0\}=\partial  M$;
\item The $1$-form $dx$ is non-vanishing on $\partial  M$.

\end{enumerate}
There exists a neighbourhood $U\subset  M$ of $\partial M$ and a
diffeomorphism $\Phi_U:U\to [0,\epsilon)\times \partial M$ such
that $x_{|U}$ is the composition of $\Phi_U$ with the
projection on the first factor. In the sequel we fix such a
product decomposition near the boundary.

The metric $g$ on $\Mci$ is said to be
of conformally cusp type if on $U\cap \Mci$ it is of the form
\begin{equation}\label{metcus}
g_p=x^{2p}\left(\frac{dx^2}{x^4}+ h\right)
\end{equation}
where $p \in (0,\infty)$ and $h$ is a metric on $\partial M$ which does not
depend on $x$. Thus $g_p= x^{2p} g_c,$ where $g_c$ is a particular
case of an exact cusp metric as in \cite{meni,wlom} and also an
exact $b$-metric in the sense of Melrose. Geometrically, $g_c$,
which takes the form $g_c = \frac{dx^2}{x^4} + h$ on $U\cap \Mci$,
is simply a metric with infinite cylindrical ends, as one can see
by switching to the variable $v=1/x$. Recall that $x$ is a global
function, thus $g_p$ is defined on $\Mci$. The motivating example
is given by complete hyperbolic manifolds of finite volume.
Outside a compact set, such a hyperbolic manifold is isometric to
an infinite cylinder $(1,\infty)\times T$ where $(T,h)$ is a
(possibly disconnected) flat manifold; the metric takes the form
\[dt^2+e^{-2t}h\]
which is easily seen to be of the form \eqref{metcus} with $p=1$ if we set $x:=e^{-t}$.

Let $E$ be a twisting bundle on $ M$, with a connection
which is flat in the direction of $\partial_x$. This
implies that near the boundary, $E$ together with its connection
are pull-backs of their restrictions to the boundary $E_{|\partial
M}$. Finally, let $D_p$ (where $2p$ is the power in the conformal
metric $g_p$) denote the twisted Dirac operator associated to the
aforementioned data.

The main assumption under which we work is the invertibility of the boundary
Dirac operator. More precisely,
\begin{assumption}\label{a1}
For each connected component $N$ of $\partial  M$, we assume that
the Dirac operator $\Dh$ on $N$ with respect to the metric $h$ and twisted by $E$,
is invertible.
\end{assumption}
Under this assumption, the results of \cite{wlom} imply that the $L^2$
spectrum of the essentially self adjoint operator $D_p$ is discrete and obeys a
Weyl-type law; moreover the eta function $\eta(D_p,s)$ is holomorphic for $\Re(s)>n$
and extends to a meromorphic function with possible \emph{double} poles
at certain points. In particular, for $n$ odd, $s=0$ is such a possible
double pole.

\begin{theorem}\label{thmm}
Under Assumption \ref{a1}, the eta function of the twisted
Dirac operator on $\Mci$ is regular at $s=0$.
\end{theorem}

It follows that we can define a ``honest" eta invariant, depending on
the eigenvalues of $D_p$, as the regular value at $s=0$ of $\eta(D_p,s)$.
\begin{proof}
We need to revisit the construction giving the meromorphic extension
and the structure of the poles of the eta function. The main tool
is the calculus of cusp pseudodifferential operators first
introduced in \cite{meni}, whose definition we review in Appendix \ref{sec-cuspcalc}. 

The spinor bundles for conformal metrics are canonically identified
together with their metrics. It follows that $D_c$, the Dirac operator
for $g_c$, is linked to $D_p$ by formula \eqref{confch}. However these two
operators act on different $L^2$ spaces because the measures $\mu_p$ and
$\mu_c$ induced by the metrics $g_p,g_c$ are not the same. We view $D_p$ as acting
in $L^2(\Mci,\Sigma\otimes E,\mu_p)$ and we conjugate it through the Hilbert space
isometry
\begin{align*}
L^2(\Mci,\Sigma\otimes E,\mu_p)\to L^2(\Mci,\Sigma\otimes E,\mu_c),&&\sigma\mapsto x^{np/2}\sigma.
\end{align*}
It follows that $D_p$ is unitarily equivalent to the operator
\[A = x^{\frac{np}{2}}D_p x^{-\frac{np}{2}} \text{ acting in
$L^2(\Mci,\Sigma\otimes E, \mu_c)$}.\]
Using the formula \eqref{confch} with $e^{-f} = x^p$, we see that
\[D_p = x^{-p\frac{n+1}{2}} D_c x^{p\frac{n-1}{2}}.\]
In the sequel we thus replace $D_p$ by the unitarily equivalent operator
\[A = x^{-\frac{p}{2}}D_c x^{-\frac{p}{2}} \text{ acting in
$L^2(\Mci,\Sigma\otimes E,\mu_c)$}.\]
This operator is an elliptic operator in the weighted cusp
calculus $x^{-p}\Diff^1_c(M, \Sigma\otimes E)$. The \emph{normal
operator} of a cusp operator $P$ in $\Diff^1_c(M, \Sigma\otimes
E)$ is the $1$-parameter family of operators on $\partial M$
defined by
\[\cN(P)(\xi)\phi=[e^{i\frac{\xi}{x}}
P(e^{-i\frac{\xi}{x}}\tilde{\phi})]_{|x=0}\]
for $\xi \in \mathbb{R}$,  where $\tilde{\phi}$ is any extension of the spinor $\phi$ from $\partial  M$
to $ M$.

Since $g_c$ and the twisting bundle $E$ and its connection are products  near infinity, we have
\begin{align}\label{dcab}
D_c=c(\nu)x^2\partial_x+\Dh,&&
A=x^{-p}\left[\left(x^2\partial_x - \frac{px}{2} \right)c(\nu)+\Dh\right]
\end{align}
where $\Dh$ is the Dirac operator on $\partial M$, $\nu=\frac{dx}{x^2}$, and $c(\nu)$ is Clifford multiplication by $\nu$. It follows from the definition that
\begin{equation}\label{cndc}
\cN(D_c)(\xi)=c(\nu) i\xi+\Dh.
\end{equation}
The boundary operator $\Dh$ anti-commutes with $c(\nu)$ for algebraic
reasons and is invertible by Assumption \ref{a1}. Since
\[
\cN(D_c)(\xi)^2 = \xi^2 + \Dh^2
\]
is strictly positive, it follows that $\cN(D_c)(\xi)$ is
invertible for all $\xi$. Such an operator is called \emph{fully
elliptic}. From \cite{wlom}, we know that $A$ has essentially the
same properties as an elliptic operator on a closed manifold: it
is Fredholm, has compact resolvent and (in the self-adjoint case)
has pure-point spectrum. The eigenvalues are distributed according
to a suitable Weyl-type law; in particular, the eta function
$\eta(A,s)$ is well-defined for large real parts of $s$. Moreover,
$A(A^2)^{-\frac{s+1}{2}}$ is a holomorphic family of cusp
operators in $x^{ps}\Psi_c^{-s}(M,\Sigma\otimes E)$ if we define
it to be $0$ on the finite-dimensional null-space of $A$, see
\cite[Proposition 15]{wlom}. It follows from \cite[Proposition
14]{wlom} that the trace of this family (i.e., the eta function)
extends meromorphically to $\cz$ with possible poles when $s \in
\{n,n-1,n-2,\ldots\}$ and when $ps \in \{1,0,-1,-2\ldots\}$; the
poles are at most double at points in the intersections of these
two sets, otherwise they are at most simple. The content of the
theorem is that $s=0$ is in fact a regular point. We will see
later that some of the above singularities do not occur in our
setting.

To start the proof, consider the holomorphic family in two complex variables
\begin{align}\label{as}
(s,w)\mapsto x^{w} A(s), && A(s)=x^{-ps}
A(A^2)^{-\frac{s+1}{2}}\in \Psi_c^{{-s}}(M,\Sigma\otimes E)
\end{align}
and the function \[F(s,w):=\Tr(x^w A(s)).\]
Clearly $F(s,ps)=\eta(D_p,s)$.

\begin{lemma} \label{lemA}
The operator $x^w A(s)$ is of trace-class for $\Re(s)>n, \Re(w)>1$.
Moreover, $F(s,w)$ is holomorphic for $\Re(s)>n, \Re(w)>1$ and
extends to $\cz\times \cz$ as a meromorphic function with possibly
simple poles in $s$ at $s\in \{n,n-1,n-2,\ldots\}$ and in $w$ at
$w\in \{1,0,-1,-2,\ldots\}$.
\end{lemma}
\begin{proof}
The operator kernel of $x^w A(s)$ is smooth outside the diagonal and
continuous at the diagonal for $\Re(s)>n$. Its restriction to the diagonal
is a smooth multiple of the cusp volume density for such $s$,
and has an asymptotic expansion in powers
of $x$ as $x\to 0$, starting from $x^w$. This is due to the fact that
$A(s)$ is a conormal distribution on the cusp double space $M^2_c$,
with Taylor expansion at the front face.

The trace of $x^wA(s)$ equals the integral on the lifted diagonal of the
above density. The normal bundle to $\Delta_c$ in $M^2_c$
is canonically identified with $\ctm$.
By the Fourier inversion formula, this is equal to
\begin{equation}\label{ft}
\int_{\ctsm} x^w a_s(p,\xi)\, \omega^n
\end{equation}
where $a_s(p,\xi)$ is a holomorphic family of classical symbols of order
$-s$ on $\ctsm$ (smooth down to $x=0$)
and $\omega$ is the canonical symplectic form on $\ctsm$.
The volume form $\omega^n$ is singular at $ M$, however $x^2$ times it
extends smoothly to the boundary of $\ctsm$. It follows that the integral
is absolutely convergent (hence holomorphic in $s,w$) for $\Re(s)>n, \Re(w)>1$.

It is now easy to construct the analytic extension of \eqref{ft} in $w$
by expanding $a_s(p,\xi)$ in  Taylor series at $x = 0$, using that
for any $k \in \nz$, we have
\[\int_0^1 x^{w + k}\, \frac{dx}{x^2} = \frac{1}{w+k-1}.\]
To get the analytic extension of \eqref{ft} in $s$ we expand $a_s(p,\xi)$
in homogeneous components in $\xi$ of order $-s - k$ where $k \in \nz$,
then switching to polar coordinates and using that
\[\int_1^\infty r^{-s - k + n - 1} dr = \frac{1}{s-n+k}.\]
\end{proof}

We first note that there is no pole at $s=0$. From here on, we assume that
the dimension of the manifold $M$ is odd, otherwise the eta function is $0$
so there is nothing to prove.
\begin{proposition}\label{p8}
The function $F(s,w)$ is regular in $s$ at $s=0$.
\end{proposition}
\begin{proof}
From the construction of the analytic extension of $F$ it follows that
for every $w$ with $\Re(w)>1$,
\[\Res_{s=0}F(s,w)=\int_{M} x^{w}
\tr(\res(A|A|^{-1})).\]
As $A$ is unitarily conjugated to $D_p$ by a real
function, we see from Corollary \ref{regconj} that the density
$\tr(\res(A|A|^{-1}))$ vanishes identically. In other words,
the holomorphic function $w\mapsto \Res_{s=0}F(s,w)$ is identically $0$
on a half-plane, and by unique continuation it is identically zero for all $w$.
\end{proof}

It remains to show that there is no pole in $w$ at $w=0$ either.
This will imply Theorem \ref{thmm} since $\eta(D_p,s)=F(s,ps)$.

For this, we fix $s$ with $\Re(s)>n$ and we examine $F(s,w)$ as a
meromorphic function in the complex variable $w$.

For any cusp operator $B\in x^z\Psi_c^{-s}(M)$ (where we suppress the bundles for brevity)
consider the power series expansion of its Schwartz kernel $x^{-z}\kappa_B$
at the front face of the cusp double space. Although
$x$ is not everywhere a defining function for the front face, we do get such an
expansion since $\kappa_B$ vanishes in Taylor series at faces other than
the front face. Taking the inverse Fourier transform of the coefficients we can regard the coefficients as lying in
the suspended calculus $\Psi^{-s}_{\rm{sus}}(\partial  M)$.
This gives a short exact sequence of spaces of operators
\begin{equation}\label{ses}
0\mapsto x^{\infty}\Psi_c^{-s}(M)\hookrightarrow x^z
\Psi_c^{-s}(M)\stackrel{q}{\to} x^z\Psi^{-s}_\sus(\partial
M)[[x]]\to 0.
\end{equation}
For a weighted cusp operator $B\in x^z \Psi_c^{-s}(M)$, we write
\[q(B)=x^z\left(q_0(B)+xq_1(B)+x^2q_2(B)+\ldots\right)\]
It is easy to see that we have $q_0(B)=\cN(x^{-z}B)$, see \cite{meni}.

We use a result from \cite{in}. Let $s>n$ and $P\in \Psi_c^{-s}(M)$ with $\Re(s)>n$. Then
$x^wP$ is trace-class for $w>1$, and for $k\in\nz$,
\begin{equation}\label{restrc}
\Res_{w=1-k} \Tr(x^wP)=\frac{1}{2\pi} \int_\rz \Tr(q_k(P)(\xi))d\xi.
\end{equation}

\begin{proposition}\label{npw} The function $F(s,w)$
does not have any poles in $w$.
\end{proposition}
\begin{proof}
Let $R\in\Diff_c^1(M,\Sigma\otimes E)$ be any cusp differential
operator which  equals $\Dp:=c(\nu)\Dh$ near the boundary. This
makes sense since we have a fixed product decomposition near the
boundary. By definition we have
\[q(R)=\Dp.\]
We notice that near $\partial M$, $\Dp$ anticommutes with the cusp differential operator $A$ from \eqref{dcab}. Therefore, if we denote by $\cI^s:=\ker(q)\subset
\Psi_c^s(M,\Sigma\otimes E)$ the subspace of operators which
vanish to every order at the front face, we have
\[RA+AR\in \cI^2.\]
This implies
\begin{align*}
[R,A^2]\in \cI^3,&& [R,(A^2)^{-\frac{s+1}{2}}]\in \cI^{-s};
\end{align*}
the latter, by the construction of the complex powers \cite{Buc}.
Together with the obvious commutation [$R,x^{-ps}]\in \cI^0,$
we get for the operator $A(s)$ defined in \eqref{as}
\[RA(s)+A(s)R\in \cI^{-s}.\]

Now for every cusp operator $Q$ we have
\begin{align*}
q(RQ)=\Dp q(Q), &&q(QR)=q(Q)\Dp
\end{align*}
because $R$ is constant in $x$ near the boundary. Therefore $\Dp
q(A(s))=-q(A(s))\Dp$. Using conjugation with the invertible
operator $\Dp$ on $\partial M$, we see that for every $\xi$,
\begin{align*}
\Tr \left(q(A(s))(\xi)\right)=\Tr\left(\Dp q(A(s))(\xi)\Dp^{-1}\right)=
-\Tr\left( q(A(s))(\xi)\Dp\Dp^{-1}\right)
\end{align*}
so $\Tr \left(q(A(s))(\xi)\right)=0$. Thus for all $k\in \nz$ the integrand
in \eqref{restrc} for $P=A(s)$ vanishes.
\end{proof}
Together with Proposition \ref{p8} this finishes the proof of
Theorem \ref{thmm} since $\eta(D_p,s)=F(s,ps)$.
\end{proof}

In fact, by invoking the regularity results of Bismut and Freed
\cite{bf}, we can restrict further the possible poles of the eta
function. By a different argument, it turns out that if $\Mci$ is
a hyperbolic manifold, then there are no poles at all (see Theorem \ref{thrw2})!
Of course, we assume that $n$ is odd since otherwise the
eta function is $0$.

\begin{theorem}\label{thrw}
Under Assumption \ref{a1}, the eta function $\eta(D_p,s)$ is
regular for $\Re(s)>-2$ and has at most simple poles at $s\in
\{-2,-4,\ldots\}$.
\end{theorem}

\begin{proof}
We have written $\eta(D_p,s)=F(s,ps)$ for an analytic function
$F(s,w)$ in $w\in\cz, s\in\cz\setminus \{n,n-1, \ldots ,
1,-1,-2,-3,\ldots \}$ by Theorem \ref{thmm}, Lemma
\ref{lemA}, and Proposition \ref{npw}.
 We claim that $F(s,w)$ is in fact regular at $s \in \cz
\setminus \{-2,-4,\ldots\}$. Indeed, for $\Re(w)>1$ the residue in
$s$ at $s = n - k$ is given by
\[ \int_M x^{w{-p(n - k)}}
\tr(\res(A|A|^{k - n -1})).\]
Since $A =
x^{\frac{np}{2}}D_p x^{-\frac{np}{2}}$, we have $A|A|^{k - n - 1} =
x^{\frac{np}{2}} D_p |D_p|^{k - n - 1} x^{-\frac{np}{2}}$, so by
\eqref{trace},
\[\Res_{s=n - k} F(s,w) = \int_M x^{w{-p(n-k)}} \tr(\res(D_p|D_p|^{k - n -1})).\]
Consider the well-known odd heat kernel small-time expansion
given in Lemma 1.9.1 of \cite{BGilkey}
\begin{equation}\label{ohk}
D_p e^{-tD_p^2}(y,y)\sim \sum_{k=0}^\infty t^{\frac{-n+k-1}{2}} b_k(y)
\end{equation}
valid on the interior of any manifold, by locality of the coefficients;
moreover, for $k$ even, we have $b_k \equiv 0$. From the relationship
\begin{equation}\label{cpho}
D_p|D_p|^{-s-1}=\frac{1}{\Gamma\left(\frac{s+1}{2}\right)}
\int_0^\infty t^{\frac{s-1}{2}} D_p e^{-tD_p^2}dt
\end{equation}
we deduce that
\begin{equation} \label{resDp}
\res(D_p|D_p|^{k - n - 1})= \frac{2}{\Gamma(\frac{n - k + 1}{2})}b_k .
\end{equation}
In particular, when $k$ is even, this residue vanishes.
By the regularity results of Bismut-Freed \cite{bf},
the pointwise traces of the
local coefficients $b_0,b_1,\ldots,b_{n}$
vanish identically, i.e. $\tr(b_k(y))\equiv 0$ for
$k=0,\ldots, n$ so for these $k$ the density
$\tr(\res(D_p|D_p|^{k - n - 1}))$ also vanishes. (The vanishing of
the term with $k=n$, corresponding to the residue of the eta function
at the origin, has already been proved in Proposition \ref{p8} by
using conformal invariance, although here we could have also deduced this fact
from the regularity results of Bismut-Freed \cite{bf}.) This proves
that $F(s,w)$ is regular at $s \in \cz \setminus \{-2,-4,\ldots\}$.
\end{proof}

In the hyperbolic case we have a stronger result. The necessary local vanishing
of the heat trace is proved in the next section.

\begin{theorem}\label{thrw2}
If $\Mci$ is an odd dimensional hyperbolic manifold of
finite volume, the eta function of the Dirac operator
twisted by a homogeneous vector bundle is entire.
\end{theorem}
\begin{proof}
By Proposition \ref{p:odd-heatkernel} in Section \ref{sec-hyper},
we have that for $m\in\mathbb{H}^n(\mathbb{R})$ with
$n=2d+1$,
\[\tr(D e^{-tD^2}(m,m))=0.\]
This implies that $\tr(b_k)=0$ for all the coefficients
$b_k$ on $\hh^n(\mathbb{R})$, and thus also on the locally
symmetric space $\Mci$. In view of \eqref{resDp}, the possible
poles in $s$ of the function $F(s,w)$ actually do not occur.
Together with Proposition \ref{npw}, this shows that the eta
function is entire since $\eta(D,s)=F(s,s)$.
\end{proof}

%%%%%%%%%%%%%%%%%%%%%%%%%%%%%%%%%%%%%%%%%%%%%%%%%%%%%%%%%%%%%%%%%%%%%%
\section{Odd heat kernels for homogeneous vector bundles over hyperbolic space}\label{sec-hyper}
%%%%%%%%%%%%%%%%%%%%%%%%%%%%%%%%%%%%%%%%%%%%%%%%%%%%%%%%%%%%%%%%%%%%%%

The real hyperbolic space $\hh^n(\mathbb{R})$ is given as the symmetric
space $\mathrm{SO}(n,1)/\mathrm{SO}(n)$. But, for our purpose, we
use the realization of $\mathbb{H}^n(\mathbb{R})$ $=$ $G/K$ where
$G=\mathrm{Spin}(n,1)$, $K=\mathrm{Spin}(n)$, which are the double
covering groups of $\mathrm{SO}(n,1)$, $\mathrm{SO}(n)$. We denote
the Lie algebras of $G$ and $K$ by $\mathfrak{g}$ and
$\mathfrak{k}$, respectively. The Cartan involution $\theta$ on
$\mathfrak{g}$ gives the decomposition
$\mathfrak{g}=\mathfrak{k}\oplus\mathfrak{p}$ where $\mathfrak{k}$
and $\mathfrak{p}$ are,
respectively, the $+1$ and $-1$ eigenspaces of $\theta$. The subspace $\mathfrak{p}$ can be identified with
the tangent space $T_o(G/K)\cong \mathfrak{g}/\mathfrak{k}$ at
$o=eK\in G/K$. Let $\mathfrak{a}$ be a fixed maximal abelian
subspace of $\mathfrak{p}$. Then the dimension of $\mathfrak{a}$
is one. Let $M=\mathrm{Spin}(n-1)$ be the centralizer of
$A=\exp(\mathfrak{a})$ in $K$ with Lie algebra $\mathfrak{m}$. We
put $\beta$ to be {the} positive restricted root of
$(\mathfrak{g},\mathfrak{a})$. Note that $A \cong \mathbb{R}$ via
$a_r=\exp(rH)$ with $H\in\mathfrak{a}$, $\beta(H)=1$.

 From now on we assume that $n$ is odd, that is, $n=2d+1$.

The spinor bundle $\Sigma$ over $\hh^n(\mathbb{R})=G/K
=\mathrm{Spin}(n,1)/\mathrm{Spin}(n)$ is defined by
\begin{equation}\label{e:spin-b}
\Sigma = \mathrm{Spin}(n,1) \times_{\tau_s} V_{\tau_s}
\longrightarrow
\hh^n(\mathbb{R})=\mathrm{Spin}(n,1)/\mathrm{Spin}(n)
\end{equation}
where $(\tau_s,V_{\tau_s})$ denotes the spin representation of
$\mathrm{Spin}(n)$. Here points of
$\mathrm{Spin}(n,1)\times_{\tau_s} V_{\tau_s}$ are given by
equivalence classes $[g,v]$ of pairs $(g,v)$ under $(gk,v)\sim
(g,\tau_s(k)v)$. In general, any $G$-homogeneous Clifford module
bundle over $\hh^n(\mathbb{R})$ is associated to
$(\tau_s\otimes\tau, V_{\tau_s}\otimes V_\tau)$ for a unitary
representation $(\tau, V_\tau)$ of $\mathrm{Spin}(n)$ as in
\eqref{e:spin-b}, which we denote by $\Sigma\otimes E$.  For
instance, the representation $\tau_s\otimes\tau_s$ of
$\mathrm{Spin}(n)$ determines a homogeneous vector bundle
$\Sigma\otimes \Sigma$ over $\hh^n(\mathbb{R})$ whose fiber is
$V_{\tau_s}\otimes V_{\tau_s} \cong \oplus_{k=0}^{d} \wedge ^k
(\mathfrak{p}\otimes\mathbb{C})$.

The space of smooth sections from $\hh^n(\mathbb{R})$ to
$\Sigma\otimes E$ is denoted by $\cun(\hh^n(\mathbb{R}),
\Sigma\otimes E)$ and can be identified with $[\cun(G)\otimes
V_{\tau_s}\otimes V_\tau]^K$ where $K$ acts on $\cun(G)$ by the
right regular representation $R$. Now a natural connection
$\nabla: \cun(\hh^n(\mathbb{R}),\Sigma\otimes E)\to
\cun(\hh^n(\mathbb{R}),\Sigma\otimes E \otimes T^*(G/K))$ is given
by
\begin{equation}\label{e:connection}
\nabla f = \sum_{i=1}^n (R(X_i)\otimes \mathrm{Id}) f\otimes X_i^*
\end{equation}
where $\{X_i\}$ is an orthonormal basis of $\mathfrak{p}$ and
$\{X_i^*\}$ is its dual basis. This connection is the unique
connection on $\cun(\Sigma\otimes E)$ which is $G$-homogeneous and
anti-commutes with the Cartan involution $\theta$ (see Lemma 3.2
of \cite{MS}). Now the Dirac operator ${D}$ on $\Sigma\otimes E$
associated to the connection $\nabla$ is defined by
\[
{D} =\sum_{i=1}^n R(X_i)\otimes c(X_i)
\]
where $c(X_i)$ denotes the Clifford multiplication.

\begin{proposition}\label{p:odd-heatkernel}
For $m\in \mathbb{H}^n(\mathbb{R})$, we have $\mathrm{tr}(D e^{-tD^2}(m,m))=0$.
\end{proposition}

\begin{proof}
Recalling $\cun(\mathbb{H}^n(\mathbb{R}),\Sigma\otimes E) \cong
[\cun(G)\otimes V_{\tau_s}\otimes V_\tau]^K$, the Schwartz kernel
of $De^{-tD^2}$ is given by a section $H_t$ in $[\cun(G)\otimes
\mathrm{End}(V_{\tau_s}\otimes V_\tau)]^{K\times K}$ satisfying
\begin{equation}\label{e:Hk}
H_t(k_1g k_2)=(\tau_s\otimes \tau)^{-1}(k_2) H_t(g)
(\tau_s\otimes\tau)(k_1)^{-1}
\end{equation}
for $k_1,k_2\in K, g\in G$, which acts on $[\cun(G)\otimes
V_{\tau_s}\otimes V_\tau]^K$ by convolution. For each $t>0$, $H_t$
lies in $[\mathcal{S}(G)\otimes \mathrm{End}(V_{\tau_s}\otimes
V_\tau)]^K$ where $\mathcal{S}(G)=\cap_{p>0}\mathcal{S}^p(G)$ with
$\mathcal{S}^p(G)$ the Harish-Chandra $L^p$-Schwartz space. For
more details, we refer to Section 3 of \cite{MS}. Taking the local
trace of $H_t$, we have that $h_t:=\tr(H_t)\in \mathcal{S}(G)$.
From \eqref{e:Hk} and recalling that a point in the homogeneous
vector bundle $\Sigma \otimes E$ is given by an equivalence class
through the relation $(gk,v)\sim(g,(\tau_s\otimes \tau)(k)v)$, we
can see that the local trace of $D e^{-tD^2}(m,m)$ is given by
$h_t(e)$ for the identity element $e\in G$.

By the Plancherel theorem (see Theorem 4.1 in \cite{Mia2}), we
have the following expression for $h_t$ at $e\in G$ (up to a
constant depending on a normalization),
\begin{equation}\label{e:Planch}
h_t(e)=\sum_{\sigma\in \hat{M}} \int^\infty_{-\infty}
\Theta_{\sigma,i\lambda}(h_t)\, p(\sigma,i\lambda) d\lambda ,
\end{equation}
where $\hat{M}$ denotes the set of equivalence classes of irreducible unitary
representations of $M$,
\[
\Theta_{\sigma,i\lambda}(h_t)=\mathrm{Tr} \int_G h_t(g)
\pi_{\sigma,i\lambda}(g)\, dg,
\]
and $p(\sigma,i\lambda)$ denotes
the Plancherel measure associated to the unitary principal
representation $\pi_{\sigma,i\lambda}$. Here the unitary principal
representation $\pi_{\sigma,i\lambda}=\mathrm{Ind}^G_{MAN}
(\sigma\otimes e^{i\lambda}\otimes\mathrm{Id})$ acts by the left
regular representation on
\[
H_{\sigma,i\lambda}=\{\ f:G\to V_\sigma\ |\
f(gma_rn)=e^{-(i\lambda+d)rH} \sigma(m)^{-1} f(g)\ \} \]
where
$n=2d+1$. By Proposition 3.6 in \cite{MS}, it follows that
\begin{equation}\label{e:Fourier}
\Theta_{\sigma,i\lambda}(h_t)= [(\tau_s \otimes \tau)|_M:\sigma](
[\sigma:\sigma_+]-[\sigma:\sigma_-])\, \lambda e^{-t\lambda^2}
\end{equation}
where $\sigma_\pm$ denotes the half spin representation of $M$
such that $\tau_s|_M=\sigma_+\oplus \sigma_-$.  By the branching
rule from $K=\mathrm{Spin}(2d+1)$ to $M=\mathrm{Spin}(2d)$ given
in Theorem 8.1.3 of \cite{Good-Wal}, we have  that for any
$\tau\in\hat{K}$, $\sigma\in\hat{M}$, $[\tau|_M:\sigma]\leq 1$ and
$[\tau|_M:\sigma]=1$ if and only if
\[
a_i-b_i\in\mathbb{Z} \quad (i,j=1,2,\ldots, d), \quad a_1\geq
b_1\geq \ldots a_{d-1}\geq b_{d-1}\geq a_d\geq |b_d|
\]
where $\tau=\sum_{i=1}^d a_ie_i$, $\sigma=\sum_{i=1}^d b_i e_i$.
Here we denote the highest weights of the representations
$\tau,\sigma$ with respect to the standard basis. This implies
that
\begin{align*}
[(\tau_s \otimes \tau)|_M:\sigma][\sigma:\sigma_+]=[(\tau_s
\otimes \tau)|_M: \sigma_+] & =[(\tau_s\otimes
\tau)|_M:\sigma_-] \\ & =[(\tau_s \otimes \tau)|_M:\sigma]
[\sigma:\sigma_-]
\end{align*} since
$\sigma_\pm=\frac12(e_1+e_2+\ldots+e_{d-1}\pm e_d)$. Now, by
Theorem 3.1 of \cite{Mia}, we also have
$p(\sigma_+,i\lambda)=p(\sigma_-,i\lambda)$. This implies
$h_t(e)=0$ by \eqref{e:Planch} and \eqref{e:Fourier}, which
completes the proof.
\end{proof}

\begin{remark} It can be proved, using the Selberg trace
formula, that the eta function vanishes at negative odd integers
under the condition $\Gamma\cap P=\Gamma\cap N$ for the
fundamental group $\Gamma$ of the given hyperbolic manifold where
$P=MAN$ (Langlands decomposition) denotes a parabolic subgroup of
$G$ fixing the infinity point of a cusp. Recall that this fact is
true for arbitrary operators of Dirac type over closed manifolds,
as follows immediately from \eqref{cpho} and from the odd heat
trace expansion \eqref{ohk}. We believe that this vanishing holds
also in the context of manifolds with conformal cusps without this
technical condition but the necessary work, which surpasses the
scope of this paper, is left for a future publication.
\end{remark}

\appendix
\section{The spectrum of the Dirac operator}\label{sec-spec}

Recall that under Assumption \ref{a1}, the Dirac operator $D_p$ is
always essentially self-adjoint with discrete spectrum
\cite{wlom}. One may ask what happens with the eta invariant when
Assumption \ref{a1} does not hold. Like in \cite{bar}, for $p\leq
1$ when the manifold is complete with respect to $g_p$, the answer
is that the continuous spectrum of the twisted Dirac operator
(which is essentially self adjoint by \cite{Wol}, \cite{St})
becomes the full real line, hence the usual definition of the eta
invariant breaks down. We will not attempt here to extend the
definition in that case, note however that for finite-volume
hyperbolic manifolds this has been done in \cite{jm}. The proof of
the following result is very similar to the corresponding
statements from \cite{MM,slf} concerning magnetic and Hodge
Laplacians.

\begin{theorem} \label{thm-spec}
Let $M$ be a spin manifold with conformal cusps and $E$ a twisting
bundle of product type near the cusps. Let $D_p$ denote the
Dirac operator associated to the metric \eqref{e:metric-p} on $M$,
twisted by $E$. If Assumption \ref{a1} does not hold, then
\begin{itemize}
\item if $0 < p\leq 1$, the essential spectrum of $D_p$ is
$\rz$.
\item if $p>1$, then $D_p$ is not essentially
self-adjoint, and every self-adjoint extension of $D_p$ in
$L^2$ has purely discrete spectrum.
\end{itemize}
\end{theorem}
\begin{proof}
The idea is to reduce the problem to a $1$-dimensional problem,
essentially to the computation of the spectrum of $i\partial_t$ on
an interval.

When $p>1$, the metric is of metric horn type, and self-adjoint
extensions of $D_p$ on $M$ (given by boundary conditions at $x=0$)
are in $1$-to-$1$ correspondence with Lagrangian subspace in
$\ker(\Dh)$ with respect to the symplectic form
\[\omega(u,v):=\langle c(\nu)u,v\rangle_{L^2(\partial M,\Sigma\otimes E)},\]
see \cite{lp}. Such subspaces exist by
the cobordism invariance of the index (note that $\partial M$ may be
disconnected). Moreover, since Assumption \ref{a1} does not hold,
there exist infinitely many Lagrangian subspaces in $\ker(\Dh)$
thus $D_p$ is not essentially self-adjoint.

We work with the operator $A$ from \eqref{dcab},
which is unitarily conjugated to ${D_p}$ hence has the same
spectrum as $D_p$. When $p>1$, for each Lagrangian
subspace $W\subset \ker(\Dh)$, $A$ is essentially self-adjoint
on the initial domain
\[\Dom_W(A)=\cun_c(M, \Sigma\otimes E)\oplus \{\phi(x) x^{p/2} w;w\in W\}\]
for some fixed cut-off function $\phi$ supported in the cusps which equals
$1$ near infinity. When $0<p\leq 1$, $A$ is essentially self-adjoint
on $\Dom(A)=\cun_c(M, \Sigma\otimes E)$.

The essential spectrum of $A$ (with the above boundary condition
when $p>1$) can be computed on the complement of
any compact set in $M$, i.e., on the union of the cusps, by
imposing self-adjoint boundary conditions. More precisely, consider
the non-compact manifold with boundary $M_\epsilon:=\{x\leq \epsilon\}$.
We need to specify a self-adjoint boundary condition for $A$ at
$x=\epsilon$, which is obtained by the $APS$ condition and by
choosing yet another Lagrangian subspace in $\ker (\Dh)$.
With these self-adjoint
boundary conditions, the decomposition principle (see \cite[Prop.\ 1]{bar})
states that the essential spectrum of $A$ on $\Me$ coincides with
the essential spectrum of $A$ on $M$.

We decompose the space of $L^2$ spinors on $\partial M$ twisted by $E$
into the space of zero-modes (i.e., the kernel of $\Dh$) and its orthogonal
complement consisting of ``high energy modes''.
Accordingly we get an orthogonal decomposition of $L^2(\Me,\Sigma\otimes E)$
into zero-modes and high energy modes, the main point being that $A$
preserves this decomposition. As in \cite[Prop.\ 5.1]{slf}
the high energy modes do not contribute to the essential spectrum.
The reason is that there exists a cusp pseudodifferential
operator $R\in x^{-2p}\Psi_c^{-\infty}(M,\Sigma\otimes E)$, localized
on the cusps and acting as $0$
on high energy modes, such that $A^2+R$ is fully elliptic. Therefore
this operator has discrete spectrum so in particular, on high energy modes
$A^2$ has discrete spectrum as claimed.
We are left with the formally self-adjoint
operator
\[
A_0:= x^{1-p}c(\nu)\left(x\partial_x -\tfrac{p}{2}\right)
\]
acting in $L^2(0,\epsilon)\otimes\ker(\Dh)$ with respect to the volume form
$\frac{dx}{x^2}$ (with a certain boundary condition at $x=\epsilon$).

We claim that $A_0$ is unitarily equivalent to $c(\nu) t\partial_t$
over a certain interval depending on $\epsilon$ and $p$,
with respect to the measure $\frac{dt}{t}$. We start by conjugating with
$x^{1/2}$, so the volume form becomes $\frac{dx}{x}$ and
\[
x^{-\frac12}A_0 x^{\frac12} = x^{2-p} c(\nu) \partial_x +
\frac{c(\nu)}{2}(1-p) x^{1-p}.
\]
For $p=1$, we already obtain the desired
expression for our operator by setting $t:=x$. For $p\neq 1$, we
write
\begin{align*}
x^{2-p}\partial_x=y^2\partial_y, && y:=(1-p)x^{1-p}.\end{align*}
Then, after conjugating with $y^{-1/2}$, we obtain the
operator $c(\nu)y^2\partial_y$ acting in $L^2$ with respect
to the measure $\frac{dy}{y^2}$. With the change of variable
\emph{cusp-to-$b$}
\[t:=e^{-1/y}\]
we get the desired operator.

Now for $p\leq 1$ the operator $c(\nu)t\partial_t$ acts on an
interval of the form $(0,\beta)$, while for $p>1$ it acts on
$(1,\beta)$ for some strictly positive $\beta$. Here
$c(\nu)$ is a diagonalizable automorphism with $\pm i$
eigenvalues. For every self-adjoint extension, in the first case
the spectrum is $\rz$ while in the second case it is discrete.
\end{proof}

Alternatively, one could prove the first part of the theorem similarly
to \cite{bar} by constructing Weyl sequences for each real number.

Although for $p>1$ the spectrum of any self-adjoint extension of
$D_p$ is discrete even when Assumption \ref{a1} does not hold, the
methods of this paper do not show the meromorphic extension of the
eta function in that case.

\section{Elements of the cusp calculus}\label{sec-cuspcalc}

\begin{figure}
\centering
\includegraphics{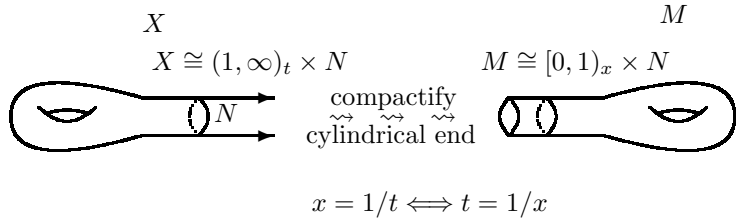}
\caption{Compactifying $X$ into $M$.} \label{fig-cusp}
\end{figure}

In this appendix we give a short introduction to the cusp calculus (first defined by Melrose and Nistor \cite{meni}). Consider a Riemannian
manifold $X$ with a cylindrical end as shown in the left-hand side
of Figure \ref{fig-cusp}.
The metric takes the form $dt^2 + h$ on the cylinder where $h$ is a metric on the cross section $N$.
Changing coordinates to $x = 1/t$ and noting that for $t \in (1,\infty)$ we have $x \in (0,1)$, and $t \to \infty$ implies $x \to 0$, it follows that we can view $X$ as the interior of the \emph{compact} manifold $M$ obtained
from $X$ by replacing the infinite cylinder $(1,\infty)_t \times N$ with
the finite cylinder $[0,1)_x \times N$, see Figure
\ref{fig-cusp}.

A cusp pseudodifferential operator is just a ``usual'' pseudodifferential operator on $X$ that ``behaves nicely'' near $t = \infty$. To make this precise, recall that in local coordinates on the cylinder, the Schwartz kernel of an $m$-th order ($m \in \mathbb{C}$) pseudodifferential operator $A$ on $X$ takes the form
\begin{equation} \label{A1}
\kappa_A = \int e^{i (t - t') \tau  + i(y - y') \cdot
\eta} \, a(t,y,\tau,\eta)\, \dbar \tau \dbar \eta ,
\end{equation}
where $(t,t',y,y') \in (1,\infty)^2 \times \mathcal{U}^2$ with $\mathcal{U}$ local coordinates on $N$ and $a$ is a classical symbol of order $m$ in $\tau$ and $\eta$.  We say that the operator $A$ is, by definition, an $m$-th order \emph{cusp
pseudodifferential operator} if in the compactified coordinates,
\[
\tilde{a}(x,y,\tau,\eta) : = a (1/x,y,\tau,\eta)
\] 
is smooth at $x = 0$ (and is still classical in $\tau$
and $\eta$). This definition only works on the cylinder, so to be more precise, $A$ is a \emph{cusp operator} if it can be written as $A_1
+ A_2 + A_3$, where $A_1$ is of the form \eqref{A1} such that
$\tilde{a}(x,y,\tau,\eta) : = a (1/x,y,\tau,\eta)$ is smooth at $x
= 0$, $A_2$ is a usual pseudodifferential operator on the compact
part of $X$, and finally, where $A_3$ is a smoothing operator on
$X$ that vanishes, with all derivatives, at $\infty$ on the
cylinder (or equivalently, the Schwartz kernel of $A_3$ is a
smooth function on $M^2$ vanishing to infinite order at $\partial
(M^2)$). The space of
cusp pseudodifferential operators of order $m\in\cz$
is denoted $\Psi_c^m$. If the symbol $a$ is polynomial in $\tau,\eta$,
the resulting operator is differential, and can be written
near $x=0$ as sums of compositions of partial
differentials on $N$ and of $\partial_t=-x^2\partial_x$
with smooth coefficients on the compactification $M$. The space of cusp differential
operators of order $m\in\zz$ is denoted ${\Diff}^m_c$.  

\begin{figure}
\centering
\includegraphics{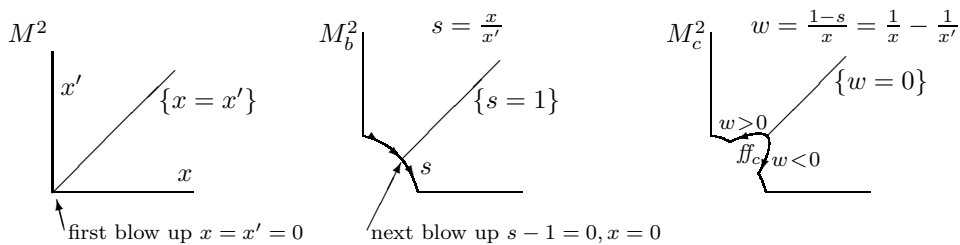}
\caption{The blown-up manifold $M^2_c$. (We omit the $N^2$ factor.)} \label{fig-M2c}
\end{figure}

Cusp operators are usually presented geometrically in relation to blown-up spaces, which might obscure their straightforward definition, so we shall explain this relationship. Setting
$z = (t - t',y - y')$, which is a \emph{normal} coordinate to the set $\{z = 0\} = \{t = t', y = y'\} =$ diagonal in $X^2$, we see from \eqref{A1} that $\kappa_A$ is written as the inverse Fourier transform of a symbol using the normal coordinate $z$. Hence $\kappa_A$ is a distribution on $X^2$ that is, by definition, \emph{conormal} to the diagonal in $X^2$. Expressing the kernel \eqref{A1} in the compactified coordinates, we obtain
\begin{equation} \label{A3}
\kappa_A = \int e^{i z \cdot (\tau, \eta)} \, \tilde{a}(x,y,\tau,\eta)\, \dbar \tau \dbar \eta \ ,\quad \text{where}\ z = \Big( \frac{1}{x} - \frac{1}{x'}, y - y' \Big).
\end{equation}
Note that $z$ is not a normal coordinate to the diagonal (given by $\{x=x',y=y'\}$) in $M^2$ because the coordinate $\frac{1}{x} - \frac{1}{x'}$ fails to be smooth at the corner $x = x' = 0$. Thus, it seems like switching to compactified coordinates destroys the conormal distribution portrayal of pseudodifferential operators. However, we now show how to interpret this kernel
as being conormal, not in $M^2$, but on a
related blown-up manifold. The idea is to blow-up the singular 
point $x = x' = 0$ until the kernel \eqref{A3} can be interpreted as
conormal. To begin this program, we first write
\[
M^2 \cong [0,1)_x \times [0,1)_{x'} \times N^2
\]
near the corner $\{x = x' = 0\}$ as shown pictorially on the left in Figure \ref{fig-M2c}.
Next, we introduce polar coordinates $(r,\theta)$ in the $x,x'$ variables, where $r = \sqrt{x^2 + (x')^2}$ and $\theta = \arctan(x'/x)$. Geometrically we can think of introducing polar coordinates as  ``blowing up'' the set $\{x = x' = 0\}$ by replacing it with a quarter-circle (the angular $\theta$ coordinate). The resulting manifold is called the \emph{$b$-double space} $M^2_b$, see the middle picture in Figure \ref{fig-M2c}. Actually, instead of using the standard polar coordinates $(r,\theta)$, in the sequel it is helpful to use the projective coordinates $(x,s)$ where $s =x/x'$, which can be used as coordinates instead of $(r,\theta)$ for
$\theta$ away from $0$ and $\pi/2$. Here, $x$ represents the radial coordinate and $s$ represents the angular coordinate along the quarter circle.

Now we form the \emph{cusp double space} $M^2_c$. To do
so, we introduce polar coordinates where $s = 1$ and $x = 0$, in $M^2_b$,
as shown in Figure \ref{fig-M2c}.
This blow-up geometrically replaces the set $s = 1,x=0$ in
$M^2_b$ with a half circle, which is called the \emph{cusp front
face} and which we denote by
$\mathit{ff}_{\hspace{-.2em} c}$.  Since the set $s = 1,x=0$ is the set of points where $1 - s = 0$ and $x = 0$, we can use the projective coordinate
\[
w = \frac{1 - s}{x} = \frac{1}{x} - \frac{1}{x'}
\]
as an angular coordinate along $\mathit{ff}_{\hspace{-.2em} c}$ and we can use $x$ as the radial variable, at least if we stay away from the extremities of $\mathit{ff}_{\hspace{-.2em} c}$.
Thus, $(x,w)$ can be used as coordinates near the blown-up face
$\mathit{ff}_{\hspace{-.2em} c}$. Note that the set $\{x =
x'\}$ corresponds to the set $\{w =
0\}$ in $M^2_c$, therefore $w$ is a
normal coordinate to $\{x = x'\}$. Moreover, in view of the
formula \eqref{A3}, we see that $\kappa_{A}$ is, by definition, a distribution
conormal to the set $\{w = 0, y = y'\}$ in $M^2_c$, which is called the \emph{cusp diagonal}. In fact, one can prove the following theorem.

\begin{theorem} The Schwartz kernels of cusp pseudodifferential operators are in one-to-one correspondence with distributions on $M^2_c$ that are conormal to the cusp diagonal and vanish to infinite order at all boundary hypersurfaces of $M^2_c$ except the cusp front face where they are smooth.
\end{theorem}

Our original definition of a cusp pseudodifferential operator as presented after \eqref{A1} is usually disregarded in favor of the more geometric definition presented in the above theorem. It is evident how to extend the definition from the scalar case
to cusp operators acting on sections of vector bundles which are
of product type in $t$ near $t=\infty$.  

\bibliographystyle{amsplain}

\end{document}